\documentclass[12pt]{article}
\usepackage{amssymb,latexsym}
\renewcommand{\thefootnote}{}

\begin{document}
\title{Bounded automorphisms and quasi-isometries of finitely generated groups}
\author{Aniruddha C. Naolekar \\  
Indian Statistical Institute, Stat-Math Unit\\ 
203 B.T.Road, Kolkata-700 108\\
and \\
Parameswaran Sankaran\\
Institute of Mathematical Sciences,\\
CIT Campus, Taramani, Chennai 600 113\\ 
Email: {\tt anirudha@isical.ac.in}\\ 
{\tt sankaran@imsc.res.in}} 
\date{}
\maketitle

\footnote{2000 A.M.S. Subject Classification:- 20F65, 20F28, 20F67\\
Key words and phrases: Automorphism, virtual 
automorphism, $FC$-centre, quasi-isometry, free product with amalgamations,
HNN extensions, hyperbolic groups.}

\thispagestyle{empty}

\def\theequation {\arabic{section}.\arabic{equation}}
\renewcommand{\thefootnote}{}

\newcommand{\codim}{\mbox{{\rm codim}$\,$}}
\newcommand{\stab}{\mbox{{\rm stab}$\,$}}
\newcommand{\lr}{\mbox{$\longrightarrow$}}

\newcommand{\ch}{{\cal H}}
\newcommand{\cf}{{\cal F}}
\newcommand{\cd}{{\cal D}}

\newcommand{\blr}{\Big \longrightarrow}
\newcommand{\da}{\Big \downarrow}
\newcommand{\ua}{\Big \uparrow}
\newcommand{\hra}{\mbox{\LARGE{$\hookrightarrow$}}}
\newcommand{\rt}{\mbox{\Large{$\rightarrowtail$}}}
\newcommand{\dua}{\begin{array}[t]{c}
\Big\uparrow \\ [-4mm]
\scriptscriptstyle \wedge \end{array}}

\newcommand{\be}{\begin{equation}}
\newcommand{\ee}{\end{equation}}
\newcommand{\ov}{\overline}
\newtheorem{guess}{Theorem}[section]
\newcommand{\bth}{\begin{guess}$\!\!\!${\bf .}~}
\newcommand{\eeth}{\end{guess}}
\renewcommand{\bar}{\overline}
\newtheorem{propo}[guess]{Proposition}
\newcommand{\bpropo}{\begin{propo}$\!\!\!${\bf .}~}
\newcommand{\epropo}{\end{propo}}

\newtheorem{lema}[guess]{Lemma}
\newcommand{\blem}{\begin{lema}$\!\!\!${\bf .}~}
\newcommand{\elem}{\end{lema}}

\newtheorem{defe}[guess]{Definition}
\newcommand{\bdefe}{\begin{defe}$\!\!\!${\bf .}~}
\newcommand{\edefe}{\end{defe}}

\newtheorem{coro}[guess]{Corollary}
\newcommand{\bcor}{\begin{coro}$\!\!\!${\bf .}~}
\newcommand{\ecor}{\end{coro}}

\newtheorem{rema}[guess]{Remark}
\newcommand{\brem}{\begin{rema}$\!\!\!${\bf .}~\rm}
\newcommand{\erem}{\end{rema}}

\newtheorem{exam}[guess]{Example}
\newcommand{\beg}{\begin{exam}$\!\!\!${\bf .}~\rm}
\newcommand{\eeg}{\end{exam}}

\newcommand{\ctext}[1]{\makebox(0,0){#1}}
\setlength{\unitlength}{0.1mm}
\newcommand{\cl}{{\cal L}}
\newcommand{\cp}{{\cal P}}
\newcommand{\ci}{{\cal I}}
\newcommand{\bz}{\mathbb{Z}}
\newcommand{\cs}{{\cal s}}  
\newcommand{\cv}{{\cal V}}
\newcommand{\ce}{{\cal E}}
\newcommand{\ck}{{\cal K}}
\newcommand{\cR}{{\cal R}}
\newcommand{\bq}{\mathbb{Q}}
\newcommand{\bt}{\mathbb{T}}
\newcommand{\bh}{\mathbb{H}}
\newcommand{\br}{\mathbb{R}}
\newcommand{\wt}{\widetilde}
\newcommand{\im}{{\rm Im}\,}
\newcommand{\bc}{\mathbb{C}}
\newcommand{\bp}{\mathbb{P}}
\newcommand{\spin}{{\rm Spin}\,}
\newcommand{\ds}{\displaystyle}
\newcommand{\tor}{{\rm Tor}\,}
\newcommand{\bs}{\mathbb{S}}
\def\ns{\mathop{\lr}}
\def\nssup{\mathop{\lr\,sup}}
\def\nsinf{\mathop{\lr\,inf}}
\renewcommand{\phi}{\varphi}
\newcommand{\co}{{\cal O}}

\noindent
{\bf Abstract:} 
Let $\Gamma$ be any finitely generated {\it infinite} group. Denote by 
$K(\Gamma)$ the $FC$-centre of $\Gamma$, i.e, 
subgroup of all elements of $\Gamma$ whose centralizers are of 
finite index in $\Gamma$. Let $QI(\Gamma)$ denote the group of 
quasi-isometries of $\Gamma$ with respect to a word metric. 
We observe that the natural homomorphism 
$\theta_\Gamma:Aut(\Gamma)\lr QI(\Gamma)$ is a monomorphism only
if $K(\Gamma)$ equals the centre $Z(\Gamma)$ of $\Gamma$. The converse 
holds if $K(\Gamma)=Z(\Gamma)$ is torsion free.  
When $K(\Gamma)$ is finite we show that $\theta_{\overline{\Gamma}}$ is 
a monomorphism where $\overline{\Gamma}=\Gamma/K(\Gamma)$. 
We apply this criterion to many interesting
classes of groups.

\section{Introduction}
Let $f:X\lr X'$ be a map (which is not assumed to be continuous)  
between metric spaces. We say that $f$ is a 
$(\lambda, \epsilon)$-quasi-isometric embedding if 
$\lambda^{-1} d(x,y)-\epsilon \leq d'(f(x),f(y)) 
\leq \lambda d(x,y)+\epsilon$ for all $x,y\in X$. Here $\lambda \geq 1, 
\epsilon\geq 0; ~d,d'$ denote the metrics on $X,X'$ respectively. 
If, further, there exists a constant $C\geq 0$ such that every 
$x'\in X'$ is within distance $C$ from the image of $f$, we say that $f$ is a 
$(\lambda,\epsilon)$-quasi-isometry. If $f$ is a 
$(\lambda,\epsilon)$-quasi-isometry then there exists a 
quasi-isometry $f':X'\lr  X$ (possibly for a different set of 
constants $\lambda',\epsilon', C'$) such that $f'\circ f $ 
(resp. $f\circ f'$)  is quasi-isometry equivalent to the identity map of 
$X$  (resp. $X'$). (Two maps $f,g: X\lr X'$ are said to be 
quasi-isometrically  
equivalent if there exists a constant $c$ such that $d'(f(x),g(x))\leq c $ 
for all $x\in X$.) Let $[f]$ denote the equivalence class of a 
quasi-isometry $f:X\lr X$. 
The set $QI(X)$ of all equivalence classes of quasi-isometries 
of $X$ is a group under composition: $[f].[g]=[f\circ g]$ for 
$[f],[g]\in QI(X)$. 
If $X'$ is quasi-isometry equivalent to $X$, then $QI(X')$ is 
isomorphic to $QI(X)$.   
     
Let $\Gamma$ be a group generated by a finite set $A$. One has the 
word metric $d_A$ (or just $d$) where $d(\gamma,\gamma')$ is the length 
(with respect to $A$) of $\gamma^{-1}\gamma'$ for $\gamma,\gamma'\in \Gamma$. 
Note that $d_A$ is left invariant. If $B$ is another finite 
generating set then the metric spaces $(\Gamma, d_A)$ and 
$(\Gamma, d_B)$ are quasi-isometric to each other.
The group $QI(\Gamma)$, which does not depend on the choice of 
the finite generating set, is an invariant of the quasi-isomorphism 
type of $\Gamma$ and is  an important object of study 
in geometric group theory initiated by M.Gromov \cite{g}. We 
refer the reader to \cite{bh} for basic facts concerning 
quasi-isometry. 

Let $f_0: \Gamma'_0\lr \Gamma_0'', f_1: \Gamma'_1\lr \Gamma''_1 $  
be isomorphisms of groups 
where $\Gamma_i', \Gamma''_i $ are 
finite index subgroups of $\Gamma$.  One has an equivalence 
relation where $f_0\sim f_1$ if there 
exists a subgroup  $\Gamma'\subset \Gamma_0'\cap \Gamma_1'$ which 
is of finite index in $\Gamma$ such that $f_0(\gamma)=f_1(\gamma) $ 
for all $\gamma\in \Gamma'$. The equivalence classes are called 
a virtual automorphisms of $\Gamma$.    
The set of all virtual 
automorphisms of $\Gamma$ form a group $Vaut(\Gamma)$ where 
$[g].[f]=[g'\circ f']$ where $f:\Gamma_0'\lr \Gamma''_0, 
g:\Gamma'_1\lr\Gamma''_1$, $f'=f\mid f^{-1}(\Gamma'')$, 
$g'=g\mid \Gamma''$ with $\Gamma''=\Gamma''_0\cap\Gamma'_1$. 
For example, it can be seen that  
$Vaut({\Bbb Z}^n)=GL(n,\bq)$. If $\Gamma$ has no finite 
index subgroups, then $Vaut(\Gamma)=Aut(\Gamma)$. 
If $\Gamma$ and $\Gamma'$ are commensurable, then 
$Vaut(\Gamma)$ and $Vaut(\Gamma')$ are isomorphic. 
For a finitely generated group $\Gamma$, it is 
easy to show that $Vaut(\Gamma)$ is countable. In general,  
these groups are expected to be `large'. However    
F.Menegazzo and J.Tomkinson 
\cite{mt} have constructed a group having 
uncountably many elements whose virtual automorphism group is trivial.   
However, no finitely generated infinite group 
with trivial virtual automorphism group seems to be known. 

Since $\Gamma'\hookrightarrow \Gamma$ is a quasi-isometry for  
any finite index subgroup $\Gamma'$ of $\Gamma$,   
any  $[f]\in Vaut(\Gamma)$  
yields an element $[f]\in QI(\Gamma)$. This 
leads to a homomorphism of groups $\eta: Vaut(\Gamma)\lr QI(\Gamma)$ 
which is easily seen to be a monomorphism. See lemma \ref{vaut} below. 
Also one has natural homomorphisms  
of groups $\sigma:Aut(\Gamma)\lr Vaut(\Gamma)$ and 
$\theta: Aut(\Gamma) \lr QI(\Gamma),$ which 
factors through $\sigma$.
It is easy to find finitely generated infinite groups
for which $\sigma$ is not a monomorphism. (For example 
let $\Gamma=\Gamma'\times \Gamma''$ where $3\leq |\Gamma'|<
\infty$ and $\Gamma''$ any finitely generated infinite group.) 
F.Menegazzo and D.J.S.Robinson \cite{mr} have 
characterized finitely generated (infinite) groups for 
which $\sigma$ is the {\it trivial} homomorphism.
We shall obtain some very general criteria for 
$\theta:Aut(\Gamma)\lr QI(\Gamma)$ to be a 
monomorphism and apply it to different classes of groups. 

Let $K(\Gamma)$ denote the the $FC$-centre of $\Gamma$, i.e., 
$K(\Gamma)$ is the subgroup of $\Gamma$ consisting of 
all $\gamma\in \Gamma$ whose 
centralizer $C(\gamma)\subset \Gamma$ 
is of finite index in $\Gamma$.  Clearly $Z(\Gamma)\subset K(\Gamma)$.

Our main results are: 

\bth \label{main}  
Let $\Gamma$ be a finitely generated infinite group. If   
$\theta:Aut(\Gamma)\lr QI(\Gamma)$ is a monomorphism
then $K(\Gamma)=Z(\Gamma)$. The converse holds if
$K(\Gamma)=Z(\Gamma)$ is torsion free. 
\eeth

\bth \label{fvc} Let $\Gamma$ be any finitely generated
infinite group such that $K(\Gamma)$ is finite.  Then
$\theta_{\overline{\Gamma}}:Aut(\overline{\Gamma})
\lr QI(\overline{\Gamma})$ is a monomorphism where
$\overline{\Gamma}=\Gamma/K$.
\eeth
 
We prove the above theorems in \S2. In \S3,  we
apply the above results to some interesting classes of
groups which arise in combinatorial and geometric group theories.

\noindent
{\bf Acknowledgments:}  The authors wish to thank the referee 
for his/her valuable comments.

\section{The $FC$-Centre}
\noindent 
We assume throughout 
that $\Gamma$ is a finitely generated {\it infinite} group 
with the word metric relative to a (fixed) finite generating set.  
If $\Gamma'\subset \Gamma$ is the inclusion of a finite index subgroup 
in $\Gamma$ we identify $QI(\Gamma')$ with $QI(\Gamma)$ via the 
isomorphism induced by the inclusion $\Gamma'\subset \Gamma$. 
 
\blem \label{vaut} Let $\Gamma$ be a finitely generated infinite group.
(i) Let $\phi:\Gamma'\lr \Gamma'' $ be an isomorphism between finite 
index subgroups $\Gamma', \Gamma''$  of $\Gamma$.  Then 
$[\phi]\in QI(\Gamma)$ is trivial if and only if the set  
$S=\{\phi(\gamma)\gamma^{-1}\mid \gamma\in \Gamma'\}$ is finite.
(ii) The natural map $Vaut(\Gamma) \lr QI(\Gamma)$ is a 
monomorphism. 
\elem
\noindent 
{\bf Proof:}(i) Note that  $S\subset \Gamma$ is finite 
$\iff$ there exists an $R>0$ such that 
$d(1,\phi(\gamma^{-1})\gamma)<R$ for all $\gamma\in \Gamma' \iff 
d(\phi(\gamma),\gamma)<R$ 
for all $\gamma\in \Gamma'$. Since $\Gamma'$ is of finite index in 
$\Gamma$, the last statement is equivalent to  
$\phi$ being  quasi-isometrically equivalent to the 
identity map of  $\Gamma$.  

\noindent
(ii) Suppose that $[\phi]=1 $ in $QI(\Gamma)$.  By part (i), the set 
$S$ is finite. We claim that $[\Gamma':H]=|S|$ where  $H=Fix(\phi)$ 
whence the element $[\phi]\in Vaut(\Gamma)$ is trivial. 
To see the claim, note that $\gamma_1H=\gamma_2H \iff  
\gamma_1^{-1}\gamma_2\in H 
\iff \phi(\gamma_1)^{-1}\phi(\gamma_2)=\gamma_1^{-1}\gamma_2 \iff 
\phi(\gamma_2) \gamma_2^{-1}=\phi(\gamma_1)\gamma_1^{-1}$.
\hfill $\Box$                 

\noindent
As an immediate consequence we get: 

\bcor \label{aut} 
The canonical homomorphism $\theta :Aut(\Gamma)\lr 
QI(\Gamma)$ is a monomorphism if 
and only if $Fix(\phi)$ has infinite index in $\Gamma$ for 
all $\phi\neq 1$ in $Aut(\Gamma)$. 
\hfill $\Box$ 
\ecor

\bdefe  (Cf. \cite{tits})  An automorphism 
$\phi:\Gamma\lr \Gamma$  is said to be {\rm bounded} 
if $Fix(\phi)$ has finite index in $\Gamma$. (Equivalently $\phi$ is bounded 
if the set $\{\phi(\gamma)\gamma^{-1}\mid \gamma\in \Gamma\}$ is finite.)  
\edefe

For $\gamma\in \Gamma$, denote by $C(\gamma)$ the centralizer of $\gamma$, 
i.e., $C(\gamma)=\{\gamma'\in \Gamma\mid~\gamma\gamma'=\gamma'\gamma\}$. 
Recall that the $FC$-centre of $\Gamma$ is the subgroup 
$K(\Gamma)=\{\gamma\in \Gamma\mid ~[\Gamma:C(\gamma)]<\infty\}$. 
The $FC$-centre is a characteristic subgroup of $\Gamma$.
It equals the subgroup of those elements of $\Gamma$ having 
only finitely many conjugates in $\Gamma$. Note that 
$\gamma\in K(\Gamma)$ if and only if conjugation by $\gamma$, $\iota_\gamma$, 
is a bounded automorphism. The set of all bounded automorphisms 
of $\Gamma$ is a normal subgroup of $Aut(\Gamma)$ which is, from 
what has been proven already, precisely $\ker(\theta)$.  

\noindent
\blem \label{central}
If an automorphism $\phi:\Gamma\lr \Gamma$ is bounded, then 
for any $\gamma\in \Gamma$, the element $\phi(\gamma)\gamma^{-1}$ 
belongs to $K(\Gamma)$. In particular, $\theta:Aut(\Gamma)\lr QI(\Gamma)$
is a monomorphism if $K(\Gamma)=\{1\}$.
\elem
\noindent
{\bf Proof:} Let $\gamma\in \Gamma$ and 
let $\gamma_1\in Fix(\phi)\cap \gamma Fix(\phi)\gamma^{-1}$.  
Then $\phi(\gamma)^{-1}\gamma_1\phi(\gamma)=\phi(\gamma^{-1})\phi(\gamma_1)
\phi(\gamma)
=\phi(\gamma^{-1}\gamma_1\gamma)=\gamma^{-1}\gamma_1\gamma$. Therefore 
$\gamma\phi(\gamma)^{-1}\gamma_1=\gamma_1\gamma\phi(\gamma)^{-1}$, 
i.e, $\gamma_1\in C(\gamma\phi(\gamma)^{-1})=C(\phi(\gamma)\gamma^{-1})$. 
Thus the subgroup $Fix(\phi)\cap \gamma^{-1}Fix(\phi)\gamma
\subset C(\phi(\gamma)\gamma^{-1})$.  It follows that  
$C(\phi(\gamma)\gamma^{-1})$ 
has finite index in $\Gamma$ and so $\phi(\gamma)\gamma^{-1}\in
K(\Gamma)$. The second assertion of the lemma now follows
immediately from corollary \ref{aut}.
\hfill $\Box$

\bcor \label{inn}   
If $\theta|Inn(\Gamma)$ is a monomorphism then $K(\Gamma)=Z(\Gamma)$. 
If, in addition,  $Z(\Gamma)$ is trivial, then $\theta$
itself is a monomorphism.
\ecor
\noindent
{\bf Proof:} Note that $\theta([\iota_\gamma])=1
\iff [\Gamma:C(\gamma)]<\infty
\iff \gamma\in K(\Gamma)$.  It follows that if $\theta|Inn(\Gamma)$ 
is a monomorphism then  $K(\Gamma)=Z(\Gamma)$.  The second 
assertion   
follows from lemma \ref{central}.
\hfill $\Box$

\noindent
{\bf Proof of Theorem \ref{main}:}
Suppose that $\theta$ is a monomorphism. The above 
corollary shows that $K(\Gamma)=Z(\Gamma)$.

Now suppose that $\Gamma$ is such that
$K(\Gamma)=Z(\Gamma)$ is torsion free.
Let $\theta([\phi])=1\in QI(\Gamma)$.
Define $\Lambda:\Gamma\lr \Gamma$ as $\Lambda(\gamma)
=\phi(\gamma)\gamma^{-1}$ for $\gamma\in \Gamma$.
By lemma \ref{central} and by our hypothesis that
$K(\Gamma)=Z(\Gamma)$, we see that $\im(\Lambda)\subset Z(\Gamma)$.
It follows that $\Lambda$ is a {\it homomorphism} of groups.
Indeed $\Lambda(\gamma_1)\Lambda(\gamma_2)
=\phi(\gamma_1)\gamma_1^{-1}\phi(\gamma_2)
\gamma_2^{-1}=\phi(\gamma_1)\phi(\gamma_2)\gamma_2^{-1}\gamma_1^{-1}
=\phi(\gamma_1\gamma_2)(\gamma_1\gamma_2)^{-1}=\Lambda(\gamma_1\gamma_2)$.
Since by lemma \ref{vaut} $\im(\Lambda)\subset Z(\Gamma)$ is a finite 
{\it subgroup} and since $Z(\Gamma)$ is torsion free by hypothesis,
it follows that $\Lambda$ is the trivial homomorphism and hence
$\phi=id$.

\noindent
{\bf Proof of Theorem \ref{fvc}:} In view of theorem \ref{main} it
suffices to show that $K(\overline{\Gamma})=\{1\}.$
To get a contradiction, assume that $K(\overline{\Gamma})\neq \{1\}$.
Then 
there exists a $\gamma\in \Gamma\setminus K$ such that $C(\gamma;K):= 
\{\gamma'\in \Gamma|[\gamma,\gamma']\in K\}$ is a finite
index subgroup of $\Gamma$.
Now let $H=C(\gamma;K)\cap(\bigcap_{1\leq i\leq k}
C(\gamma_i))$, where $K=\{\gamma_1,\cdots, \gamma_k\}$. 
Note that $H$ is of finite index in 
$\Gamma$.  Now the map $\gamma'\mapsto [\gamma,\gamma']$ defines  
a {\it homomorphism} $\psi: H \lr K$. Since $K$ is finite, $Ker(\psi)$ 
has finite index in $H$ and hence in $\Gamma$. But then 
$Ker(\psi)$ is evidently a subgroup of  
$C(\gamma)$. This forces $\gamma\in K$, contrary to our hypothesis. 
Hence $K(\overline{\Gamma})=1$. By lemma \ref{central}, 
$\theta:Aut(\bar{\Gamma})\lr QI(\bar{\Gamma})$ is a monomorphism. 
\hfill $\Box$

\section{Examples}
There are many interesting class of  
finitely generated infinite groups for which 
$\theta_\Gamma:Aut(\Gamma)\lr QI(\Gamma)$ is a monomorphism.
In this section we give examples which arise in 
combinatorial group theory and geometric topology. Throughout 
we assume that $\Gamma$ is finitely generated and infinite.

\beg {\bf Hyperbolic groups, $C'(1/6)$-groups, etc.} 
Clearly if $\Gamma$ is a (finitely generated) group having no 
finite quotients then $\theta_\Gamma$ is a monomorphism. 
B.Trufflaut \cite{tr} has shown that in any $C'(1/6)$-group 
$\Gamma$, centralizer of any non-trivial element is cyclic. 
It follows that if $\Gamma$ is 
such a group and is not virtually cyclic then $K(\Gamma)=1$ 
and so $\theta$ is a monomorphism. 
Note that if $\Gamma$ does not contain a subgroup 
isomorphic to $\bz^2$ and is not virtually cyclic (e.g. 
a non-elementary (word) hyperbolic group) then $K(\Gamma)$ 
is a torsion group. If $K(\Gamma)$ is finitely generated, 
then it has to be finite since every element of $K(\Gamma)$ has 
only finitely many distinct conjugates. It follows 
that $\theta_{\ov{\Gamma}}$ is a monomorphism. 

Now suppose $\Gamma$ is $\delta$-hyperbolic. 
We claim that $K(\Gamma)$ is finite, 
so that $\theta_{\ov{\Gamma}}$ is a monomorphism.
Indeed, any element of $K(\Gamma)$,  being torsion, is conjugate 
to an element in the  $(4\delta +2)$-ball $B$ about $1\in \Gamma$.
(See p. 460, \cite{bh}.) Since $K(\Gamma)$ is normal, it follows that 
$K(\Gamma)$ is the union of $\Gamma$-conjugates of the {\it finite} set  
$B\cap K(\Gamma)$. Since each element of 
$K(\Gamma)$ has only finitely many conjugates in $\Gamma$, it 
follows that $K(\Gamma)$ has to be finite. 
(When $\Gamma$ is torsion free,  the assertion that $\theta_\Gamma$ 
is a monomorphism also follows from theorem \ref{main} 
since $K(\Gamma)=1=Z(\Gamma)$.) 
\eeg

However, if $\Gamma$ 
is a torsion free CAT(0) group, then $\theta_\Gamma$ 
is not always a monomorphism. For example when $\Gamma$ 
is the fundamental group of the Klein bottle, $\theta_\Gamma$ 
is not a monomorphism. See remark \ref{klein} below. 

\beg {\bf Lattices in Lie groups} \label{lie}
(i) Let $\Gamma$ be a lattice in a simply connected nilpotent 
Lie group $N$.  We claim that $\theta $ is a monomorphism. 
Indeed if $\phi\in Aut(\Gamma)$ fixes a finite index 
subgroup $\Gamma'$ of $\Gamma$, then $\Gamma'$ is also 
a lattice in $N$.  It follows from Cor. 1, p. 34 of \cite{rag} that 
the unique extension $\wt{\phi}:N\lr N$ of $\phi$ has to be 
the identity of $N$. Thus $\phi$ itself has to be 
the identity automorphism. Therefore $\theta:Aut(\Gamma)\lr 
QI(\Gamma)$ is a monomorphism.

(ii)  Let $\Gamma$ be a lattice in a real semisimple 
Lie group $G$ without compact factors.    
Then the $FC$-centre of $\Gamma$ is contained in the 
centre of $G$ by Cor. 5.18, p.84, \cite{rag}.  It follows that 
$K(\Gamma)=Z(\Gamma)\subset Z(G)$. Since the centre of $G$ is 
a finite subgroup, theorem \ref{fvc} implies that 
$\theta_{\ov{\Gamma}}$ is a monomorphism.  

When 
$\Gamma\subset G$ is an irreducible non-uniform
lattice and $G\neq SL(2,\br)$, it is known from the work of B.Farb,
R.Schwartz, and A.Eskin that
$QI(\Gamma)$ is isomorphic to the commensurator of $\Gamma$ in $G$,
which is actually {\it equal} to $Vaut(\Gamma)$. See \cite{farb} 
and the references therein. This is due to
the fact that any virtual automorphism of $\Gamma$ can be lifted 
to an automorphism of $G$. 
When $G=SL(2,\br)$, then $\Gamma$ has a finite index subgroup
which is free of finite rank is hence quasi-isomorphic to
a free group of rank $2$. (When $\Gamma$ is a uniform lattice,
then by Svarc-Milnor lemma, $QI(\Gamma)=QI(G)=QI(G/K)$ where
$K$ is a maximal compact subgroup of $G$.)
\eeg

\beg  {\bf Infinite dihedral group}\label{dihedral}
When $|\Gamma_i|=2$ for $i=1,2$,  
$\Gamma$ is the infinite dihedral group $\langle x,y|x^2=1, 
xyx=y^{-1}\rangle$. The $FC$-centre of $\Gamma$ is the 
infinite cyclic group generated by $y$ whereas the centre is 
trivial. Thus $\theta$ is not a monomorphism. 
Indeed the kernel of $\theta$ is the infinite cyclic group  
generated by the automorphism $\phi:\Gamma\lr \Gamma$ defined as 
$\phi(x)=xy, \phi(y)=y$. 
\eeg

Note that the infinite dihedral group is 
isomorphic to the free product $(\bz/2\bz)*(\bz/2\bz)$.

\beg {\bf Free product with amalgamation} 
Suppose that $\Gamma$ decomposes non-trivially as a free product with 
amalgamation: $\Gamma_1*_{\Gamma_0}\Gamma_2$ where at least
one of the indices 
$|\Gamma_i:\Gamma_0|\geq 3.$ Then\\
\noindent 
{\bf Claim:} 
$K(\Gamma)=Z(\Gamma)=Z(\Gamma_1)\cap Z(\Gamma_2)\cap \Gamma_0$.\\
\noindent
{\bf Proof:} 
We begin by showing that $K(\Gamma)\subset \Gamma_0$. 
Assume that $\gamma$ is cyclically reduced.  If $\gamma\notin \Gamma_1
\cup \Gamma_2$, we can write $\gamma=\gamma_1\cdots \gamma_{2k}$ where 
no two successive elements in the expression 
belong to the same $\Gamma_j$.  
Replacing $\gamma$ by $\gamma^{-1}$ if necessary, we have 
$\gamma_{2i-1}\in\Gamma_1\setminus \Gamma_0, \gamma_{2i}\in \Gamma_2\setminus 
\Gamma_0.$     
Since $[\Gamma_1:\Gamma_0]\geq 3$, there exists an $a\in \Gamma_1$ 
such that $a^{-1}\gamma_1\notin \Gamma_0$.  Let $b\in 
\Gamma_2\setminus\Gamma_0$.  Then for any $m>0$, $(ab)^{-m}\gamma (ab)^m$ 
is reduced as written. So $\gamma$ has infinitely many conjugates. 
This is a contradiction.  
  
Assume $\gamma=\gamma_1\in \Gamma_1\setminus \Gamma_0$.
With $a$, $b$ as in the previous paragraph and $m$ any positive integer, 
$(ba)^{-m}\gamma(ba)^m$ is reduced as written so we have a 
contradiction again.  Similarly if $\gamma\in \Gamma_2\setminus \Gamma_0$ 
we see that $(ab)^{-m}\gamma(ab)^m$ is a reduced expression 
as written for all positive integers $m$, again a contradiction. 

To complete the proof of the claim we must show $\gamma\in Z(\Gamma)$.
To get a contradiction, assume that $x_i\gamma=\gamma x'_i$, 
$x_i,x'_i\in \Gamma_i\setminus \Gamma_0, i=1,2$   
and $x_j\neq x'_j$ for 
some $j$. Then $\gamma^{-1}x_1x_2\gamma 
=x_1'x_2'$ and hence 
$\gamma^{-1}(x_1x_2)^n\gamma=(x_1'x_2')^n\neq (x_1x_2)^n$ for any 
positive integer $n>0$.  
This implies that $\gamma\notin K(\Gamma)$, contrary to our 
assumption. Hence $\gamma\in Z(\Gamma)$. 
The equality $Z(\Gamma)=Z(\Gamma_1)\cap Z(\Gamma_2)\cap \Gamma_0$   
holds since $\Gamma$ is generated by $\Gamma_1\cup \Gamma_2$. \hfill $\Box$ 

It follows from theorem \ref{main} that 
if $Z(\Gamma_1)\cap Z(\Gamma_2)\cap \Gamma_0$ is  
torsion free then $\theta$ is a monomorphism. 
\eeg 

Let $\Gamma$ be an HNN extension $\Gamma_1*_{\Gamma_0}
=\langle \Gamma_1, t\mid t^{-1}\gamma t=\phi(\gamma), \gamma\in 
\Gamma_0\rangle$.
If $\Gamma_0=\phi(\Gamma_0)=\Gamma_1$ then $\Gamma$ is 
isomorphic to a semi-direct product of $\Gamma_1$ with the 
infinite cyclic group $\langle t
\rangle$. If one of $\Gamma_0, \phi(\Gamma_0)$ -- say $\Gamma_0$ 
is a proper subgroup of $\Gamma_1$, then one has the 
following description of the centre of $\Gamma$.

\blem 
With above notations, assume that $[\Gamma_1:\Gamma_0]>1$. 
Then $Z(\Gamma)=Z(\Gamma_1)\cap Fix(\phi)$.
\elem 
\noindent 
{\bf Proof:} It is easily seen that $Z(\Gamma_1)\cap Fix(\phi)\subset 
Z(\Gamma).$ Also it is clear from our hypotheses that $t^p\notin Z(\Gamma)$
for any $p\neq 0$. 
Let $\gamma\in Z(\Gamma)$. If $\gamma\in \Gamma_1$,  then it is easy to 
show that $\gamma\in Z(\Gamma_1)\cap Fix(\phi)$. 
Choose right coset representatives for 
$\Gamma_0\backslash\Gamma_1$ and $\phi(\Gamma_0)\backslash\Gamma_1$. 
We choose $1\in \Gamma_0$ as the representative for both $\Gamma_0$ 
and $\phi(\Gamma_0)$. Write $\gamma$ in normal form 
$\gamma=\gamma_0t^{n_1}\gamma_1\cdots t^{n_k}\gamma_k $ where 
$\gamma_0\in \Gamma_1$, 
$\gamma_i,i\geq 1$ are among the chosen coset representatives, and  
$n_k\neq 0,$ for $k\geq 1$. (See p. 181 \cite{ls}.) 
Assume that $\gamma\notin \Gamma_1$ so that 
$k\geq 1$. Since $\gamma=\gamma_k\gamma\gamma_k^{-1}$, 
by uniqueness of the normal form we conclude that $\gamma_k=1$. 

We claim that $\gamma_0\in Z(\Gamma_1)\cap Fix(\phi)$.
Indeed, write  
$\gamma_0=\gamma'\gamma_0'$ for some $\gamma'\in \Gamma_0$ and $\gamma_0'$ 
chosen coset representative of $\Gamma_0\gamma_0$. Then 
$\gamma=t^{-1}\gamma t$ has normal form 
$\phi(\gamma')t^{-1}\gamma_0't^{n_1}\gamma_1\cdots t^{n_k+1}$ 
which again contradicts the uniqueness of the normal form unless 
$\gamma_0'=1, \phi(\gamma')=\phi(\gamma_0)=\gamma_0$ and $k=1$.  Thus 
$\gamma_0\in Fix(\phi)$. Thus $\gamma_0t=t\gamma_0$. In particular, 
$\gamma^{-1}$ has normal expression $\gamma^{-1}=\gamma_0^{-1}t^{-n_1}$. 
We may therefore assume that $n_1<0$. 
Since $\gamma_0\notin Z(\Gamma_0)\cap Fix(\phi)$ and $\gamma_0\in 
Fix(\phi)$, we conclude that $\gamma_0\notin Z(\Gamma_1)$. 
Choose a $\Gamma_0$-coset 
representative $\gamma''$ which does not commute with $\gamma_0\in \Gamma_0$. 
Now the normal expression for $\gamma=\gamma''^{-1}\gamma\gamma''$ is 
$ \gamma''^{-1}\gamma_0t^{n_1}\gamma''$. This again contradicts 
the uniqueness of the normal form. This shows that 
$\gamma_0=1$

Now $\gamma=t^{n_1}\gamma_1\cdots t^{n_k}$, and  $k\geq 2$ 
as $t^p\notin Z(\Gamma)$ for $p\neq 0$. 
It is now obvious that $\gamma\neq t^{-1}\gamma t$ as the normal form 
for $t^{-1}\gamma t$ is $t^{n_1-1}\gamma_1\cdots t^{n_k+1}$ 
which is evidently different from that of $\gamma$. 
This contradiction shows that $\gamma\notin \Gamma_1$. 
\hfill $\Box$

\brem \label{klein}
In case $\Gamma_1=\Gamma_0=\phi(\Gamma_0)$, then the above 
description of the centre is not valid. For example 
if $\Gamma=\bz*_\bz$ is the fundamental group of the 
Klein bottle, namely, 
$\langle t,\gamma\mid t^{-1}\gamma t
=\gamma^{-1}\rangle$, 
then $Fix(\phi)=1$, whereas $Z(\Gamma)=\langle 
t^2\rangle$. Note that $K(\Gamma)=\langle t^2,\gamma\rangle
\neq Z(\Gamma)$. 
\erem

We shall now proceed to describe the $FC$-centre of $\Gamma$.
For this purpose, we define certain subgroups 
$\Gamma_{0,k}\subset 
\Gamma_0$ on which one has a well-defined monomorphism $\phi^k$ 
as follows: $\Gamma_{0,1}=\Gamma_0$, 
and, for $k\geq 2$, 
$\Gamma_{0,k}:=\phi^{-1}(\Gamma_{0,k-1})\subset \Gamma_0$.
Set $F_k:=\{\gamma\in \Gamma_{0,k}\mid \phi^k(\gamma)=\gamma\}$ 
and let $F_\infty =\bigcup_{k\geq 1} F_k.$  Note that $F_\infty$ 
is a subgroup of $\Gamma_0$. Using Britton's lemma it is 
easy to show that for any $\gamma\in \Gamma$, one has  
$\gamma\in F_k$ if and only if $t^k\in C(\gamma)$. 

\blem 
With notation as above, assume that $[\Gamma_1:\Gamma_0]>1$. 
Then \\ (i) $K(\Gamma)\subset K(\Gamma_1)\cap F_\infty.$\\
(ii) If $F_\infty=Fix(\phi)$ is a normal subgroup of $\Gamma_1$, 
then $K(\Gamma)=K(\Gamma_1)\cap Fix(\phi)$.
\elem 
\noindent
{\bf Proof:} (i)  Suppose $\gamma\in K(\Gamma)$. Then $\gamma 
t^{k}\gamma^{-1}t^{-k}=1$ for some integer $k\neq 0$.  
Using Britton's lemma it is easy to see 
that $\gamma\in F_k\subset F_\infty$. 
As it is obvious that $K(\Gamma)\cap \Gamma_1\subset K(\Gamma_1)$ 
it follows that $K(\Gamma)\subset K(\Gamma_1)\cap F_\infty$.    

\noindent
To prove (ii), suppose $F_\infty=Fix(\phi)$ and let 
$\gamma\in K(\Gamma_1)
\cap Fix(\phi)$. Since $Fix(\phi)$ and $K(\Gamma_1)$ are 
normal in $\Gamma_1$, 
all $\Gamma_1$-conjugates of 
$\gamma$ are in $K(\Gamma_1\cap Fix(\phi)$.   
In particular, every $\Gamma_1$ conjugate of $\gamma$ commutes  
with $t$. Since $\Gamma$ is generated by 
$\Gamma_1$ and $t$, it follows that the set of all 
$\Gamma$-conjugates 
of $\gamma$ is the {\it same} as the set of all 
$\Gamma_1$-conjugates. As $\gamma\in K(\Gamma_1)$, 
latter set being finite, we conclude that $\gamma\in K(\Gamma)$.
\hfill $\Box$


\end{document}